\documentclass[a4paper,11pt,english]{amsart}
\usepackage{amssymb,babel,xspace,amscd}
\usepackage[matrix,arrow,curve]{xy}\CompileMatrices
\usepackage{hyperref}

\title[Homogeneous CR manifolds and CR algebras]{On homogeneous CR manifolds\\ and their CR algebras}

\author{Andrea Altomani}
\address{A.\ Altomani:
Dipartimento di Matematica\\ II Universit\`a di Roma
``Tor Ver\-ga\-ta''\\ Via della Ricerca Scientifica\\ 00133 Roma
(Italy)}
\email{altomani@mat.uniroma2.it}

\author{Costantino Medori}
\address{C.\ Medori:
Dipartimento di Matematica\\ Universit\`a di Parma\\ Parco Area
del\-le Scien\-ze 53/A \\ 43100 Parma (Italy)}
\email{costantino.medori@unipr.it}

\date{28 October 2005; Revised version: 18 December 2005.}

\numberwithin{equation}{section}

\theoremstyle{plain}
\newtheorem{thm}{Theorem}[section]
\newtheorem{lem}[thm]{Lemma}

\newtheorem{prop}[thm]{Proposition}

\theoremstyle{definition}
\newtheorem{exam}[thm]{Example}

\newcommand{\CR}{\texorpdfstring{\ensuremath{CR}}{CR}\xspace}
\newcommand{\C}{\mathbb{C}}
\newcommand{\R}{\mathbb{R}}

\newcommand{\gr}{\mathbf}
\newcommand{\G}{\gr{G}}

\newcommand{\la}{\mathfrak}
\newcommand{\g}{\la{g}}
\newcommand{\q}{\la{q}}
\newcommand{\f}{{\la{f}}}
\newcommand{\h}{\la{h}}

\newcommand{\cx}[1]{\Hat{#1}}

\newcommand{\ro}{\mathcal{R}}
\newcommand{\B}{\mathcal{B}}
\newcommand{\Q}{\mathcal{Q}}

\newcommand{\X}{{\mathbf X}}
\newcommand{\HHH}{{\mathcal H}}

\DeclareMathOperator{\supp}{supp}

\begin{document}

\begin{abstract}
In this paper we show some results on homogeneous \CR manifolds,
proved by introducing their associated \CR algebras.
In particular, we give different notions
of nondegeneracy (generalizing the usual notion for the Levi form)
which correspond to geometrical properties for the corresponding manifolds.
We also give distinguished equivariant \CR fibrations for homogeneous
\CR manifolds.
In the second part of the paper we apply these results to minimal
orbits for the action of a real form of a semisimple Lie group $\hat\G$
on a flag manifold $\hat\G/\mathbf Q$.
\end{abstract}

\maketitle


\section*{Introduction}
The study of the automorphisms group of \CR manifolds
was pioneered in 1932 by \'E. Cartan, who investigated the
holomorphic invariants of real hypersufaces in $\C^2$.
It was carried on
for real hypersurfaces in $\C^n$
by Morimoto and Nagano (see \cite{MoNa:1963}) and
continued by Chern \& Moser and N.Tanaka in the seventies (see
in particular \cite{ChMo:1974} and \cite{Ta:1976}).
\par
In this context homogeneous \CR manifolds are natural objects to
consider.
Orbits for the action of a real Lie group of biolomorphisms of a
complex manifold provide large classes of examples.
In particular, the orbits of a real form $\G$ of a complex semisimple
Lie group $\hat\G$ in a flag manifold $\hat\G/\mathbf Q$ have
been
investigated by Wolf \cite{Wolf:1969}
and later by several other authors, also
in connection with representation theory (see for example \cite{BKZ:1992}, \cite{HW:2003}).
Only one orbit of $\G$ in $\hat\G/\mathbf Q$  is compact.
It has minimal dimension and is called \emph{minimal orbit}.
\par
In \cite{Ta:1970}
N.Tanaka associated to any (regular) \CR manifold $M$ with
non-degenerate Levi form
a graded Lie algebra $\g$ and,
by a generalization of the theory of prolongation (in the sense of Sternberg),
proved that the group of \CR automorphisms of $M$
is a Lie group whose dimension is bounded by $\dim\g$
(see also \cite{AS:2001}).
Maximally homogeneous \CR manifolds
were studied by Nacinovich and the second author
(see \cite{MN:1997},
\cite{MN:2000} and \cite{MN:2001}).
They are compact if and only if are minimal orbits, but not all minimal orbits
are maximally homogeneous.
The notion of \CR algebra
and more general conditions of nondegeneracy were introduced in
\cite{MN:2005},
allowing an in-depth study of the \CR structure of minimal orbits
(see \cite{AMN:2005}).
\par
On the other hand, the study of compact homogeneous \CR manifolds of hypersurface type
was carried on by Azad, Huckleberry \& Richtofer \cite{AHR:1985}
and, more recently, by Alekseevsky \& Spiro (see \cite{Sp:2000}, \cite{AS:2002a} and \cite{AS:2003}), who classified all the simply connected ones.
\par\medskip
In this paper we present some recent results
on homogeneous \CR manifolds
obtained by Nacinovich and the authors
by introducing \CR algebras (see \cite{MN:2005, AMN:2005}).
\par
In the first part (\S\ref{secprelim}--\ref{secequifib}) we deal
with general homogeneous \CR manifolds. After setting basic
definitions and notation, we introduce
the \emph{\CR algebra} associated to a homogeneous \CR manifold $M$
(\S\ref{sec:homCR}). This is a pair $(\g,\q)$ consisting of a
real Lie algebra $\g$ and a subalgebra $\q$ of its
complexification: it gives an infinitesimal description of the
\CR structure of $M$.
\par
Next (\S\ref{sec:nondegeneracy}) we recall the definition of a \CR
manifold of finite type and three different nondegeneracy
conditions, to introduce analogous definitions for \CR algebras
and state an equivalence result between the conditions on \CR
manifolds and those on \CR algebras.
\par
In \S\ref{secequifib} we introduce the important notion of equivariant
\CR fibrations for manifolds and algebras. In particular we describe
the fundamental and the weakly nondegenerate reductions of \CR
algebras. For (locally) homogeneous \CR
manifolds they correspond to local equivariant \CR fibrations
that generalize the classical
Levi foliation of Levi-flat \CR manifolds.
\par\smallskip
The second part of the paper
(\S\ref{sec:flag}--\ref{sec:fundamental}) consists of applications
of these techniques to  \emph{minimal orbits} of a real semisimple
Lie group $\G$ in a flag manifold $\hat\G/\mathbf Q$ of its
complexification.
\par
The corresponding \CR algebras (\emph{parabolic minimal}) can be
effectively described in terms of the root system of $\hat\g$ with
respect to a suitably chosen Cartan subalgebra, and are classified by
their \emph{cross marked Satake diagrams}
(\S\ref{sec:parabolic}--\ref{secsatake}).
\par
In the last two sections
(\S\ref{sec:crfibr}--\ref{sec:fundamental}) we describe equivariant \CR
fibrations of minimal orbits with the aid of cross marked Satake diagram
and show that the fiber is still a minimal orbit, whose Satake diagram
can be computed from those of the base and the total space.
Then we focus on the fundamental reduction.
\section{Preliminaries on \CR manifolds}\label{secprelim}
An (abstract) almost \CR manifold of type $(n,k)$ is a triple
$(M,HM,J)$, consisting of a paracompact smooth manifold $M$ of
real dimension $(2n+k)$, a smooth subbundle $HM$ of $TM$ of even
rank $2n$ (its \emph{holomorphic tangent space}) and a smooth
partial complex structure $J:HM\to HM$, $J^2=-\mathrm{Id}$, on the
fibers of $HM$. The integer $n\geq 0$ is the \emph{\CR dimension}
and $k$ is the \emph{\CR codimension} of $(M,HM,J)$.
\par
Let $T^{1,0}M$ and
$T^{0,1}M$ be the complex subbundles of the complexification $\C\otimes HM$
of $HM$, which correspond to the $i$- and
$(-i)$-eigenspaces of $J$:
\[
T^{1,0}M\,=\,\{X\,-\, i\,J\,X\, \big{|}\, X\in HM\}
\,,\quad
 T^{0,1}M\,=\,\{X\,+\, i\,J\,X\, \big{|}\, X\in HM\}\,.
\]
We say
that $(M,HM,J)$ is a \emph{\CR manifold} if the formal integrability
condition:
\begin{equation}\label{eq:intcondit}
\left[ \mathcal C^\infty(M,T^{0,1}M),\,
\mathcal C^\infty(M,T^{0,1}M)\right]\,
\subset\, \mathcal C^\infty(M,T^{0,1}M)
\end{equation}
{holds} (we get an equivalent condition by
substituting $T^{1,0}$ for $T^{0,1}$ in \eqref{eq:intcondit}).
When $k=0$, we have $HM=TM$ and,
via the Newlander-Nirenberg theorem, we recover the
definition of a complex manifold. A smooth real manifold
of real dimension $k$ can
always be considered as a \emph{totally real} \CR manifold, i.e.
a \CR manifold of
\CR dimension $0$ and \CR codimension $k$.
\par
\smallskip
Let $(M_1,HM_1,J_1)$, $(M_2,HM_2,J_2)$ be two abstract
smooth \CR manifolds.
A smooth map $f:M_1\to M_2$,
with differential $f_*:TM_1\to TM_2$,
is a \emph{\CR map} if $f_*(HM_1)\subset HM_2$,
and $f_*(J_1v)=J_2f_*(v)$ for every $v\in HM_1$.
We say that $f$ is a \emph{\CR diffeomorphism} if
$f:M_1\to M_2$ is a smooth diffeomorphism and
both $f$ and $f^{-1}$ are \CR maps.\par
A \emph{\CR function} is a
\CR map $f:M\to\C$
of a \CR manifold $(M,HM,J)$ to $\C$, endowed with the standard
complex structure.
\par
A \CR \emph{embedding} $\phi$ of an abstract \CR manifold
$(M,HM,J)$ into a complex manifold $\X$, with complex structure
$J_{\X}$, is a \CR map which is a smooth embedding and satisfies
$\phi_*(H_pM)=\phi_*(T_pM)\cap J_{\X}(\phi_*(T_pM))$ for every
$p\in M$.
\par
If $\phi:M\to\X$ is a smooth embedding of a paracompact smooth
manifold $M$ into a complex manifold $\X$, for each point $p\in M$
we can define $H_pM$ to be the set of tangent vectors $v\in T_pM$
such that $J_{\X}\phi_*(v)\in\phi_*(T_pM)$. For $v\in H_pM$, let
$J_Mv$ be the unique tangent vector $w\in H_pM$ satisfying
$\phi_*(w)=J_{\X}\phi_*(v)$. If the dimension of $H_pM$ is
constant for $p\in M$, then $HM={\bigcup}_{p\in M}{H_pM}$ and
$J_M$ are smooth and define the unique \CR structure on $M$ for
which $(M,HM,J_M)$ is a \CR manifold and $\phi:M\to \X$ is a \CR
embedding. In particular, this is the case when $M$ is an orbit
for the action of a Lie group of biholomorphisms of $\X$.
\par
In the next sections, to shorten notation, we shall write simply
$M$, or $M^{n,k}$, for a \CR manifold $(M,HM,J)$ of type $(n,k)$,
when the \CR structure will be clear from the context.
Moreover, we shall denote by $\hat{V}$ and $\hat{\phi}$ the
complexification of a vector space $V$ and a linear map $\phi$,
respectively.
\section{The \CR algebra associated to a homogeneous \CR manifold}
\label{sec:homCR}
Let $(M,HM,J)$ be a \CR manifold, which is homogeneous for the action
of a real Lie group $\G$ of \CR transformations, so that $M=\G/\G_+$
for a closed isotropy subgroup $\G_+$  of $\G$. Set $eG_+=o\in M$.
\par
To $M$ we can associate a
\emph{\CR algebra} $(\g,\q)$.
This is a pair consisting of the real Lie algebra $\g$ of
the group $\G$ and of a complex subalgebra $\q$ of its
complexification $\hat\g$.
This subalgebra is the inverse image
\begin{equation}
\q=\cx{\pi}_*^{-1}(T_o^{0,1}M)
\end{equation}
by the complexification
$\cx{\pi}_*$ of the differential $\pi_*:\g\simeq
T_e\G \to T_oM$ of the group action at $e$.
The fact that $\q$ is a complex Lie subalgebra of
$\hat\g$ is a consequence of the formal integrability condition
\eqref{eq:intcondit}
for \CR manifolds.
\par
The real Lie subalgebra $\g_+=\g\cap\q$ is called \emph{isotropy}
and the subspace $\HHH_+=(\q+\bar\q)\cap\g$ is called
\emph{holomorphic tangent space} (where the conjugation is taken
with respect to the real form $\g$ of $\hat\g$).
\par
Vice versa, let $(\g,\q)$ be a \CR algebra.
Let $\G$ be a connected and simply connected Lie group with
Lie algebra $\g$.
Assume that the analytic subgroup $\G_+$
of $\G$ corresponding to the Lie subalgebra $\g_+=\g\cap\q$
is closed in $\G$.
Then the homogeneous space $M=\G/\G_+$ is
a smooth paracompact manifold and
has a unique \CR structure such that\,:
\begin{itemize}
\item $T^{0,1}_oM=\cx{\pi}_*(\q)$;
\item $\G$ acts on $M$ by \CR diffeomorphisms.
\end{itemize}
We denote the \CR manifold
$\G/\G_+$
by $\tilde{M}(\g,\q)$.
\par\smallskip
Let $(\g,\q)$ be a \CR algebra of type $(n,k)$.
We call $(\g,\q)$:
\begin{itemize}
\item
\emph{totally real} if
$\q=\hat\g_+=\C\otimes_{\R}\g_+$, i.e. if
$\g_+={\HHH }_+$ or, equivalently, the \CR dimension~$n$ is~$0$;
\item
\emph{totally complex} if $\q+\overline{\q}=\hat\g$,
i.e. if
${\HHH }_+=\g$, or,
equivalently, the \CR codimension~$k$ is~$0$;
\item
\emph{transitive}, or \emph{effective},
if $\g_+$ does not contain any nonzero ideal of
$\g$.
\end{itemize}
If $(\g,\q)$ is the $CR$ algebra associated to a $\G$-homogeneous
$CR$ manifold $M$, these notions express geometric properties of
$M$ (see \cite{MN:2005}). In particular, effectiveness is
equivalent to almost effectiveness of the $\G$ action (i.e.\
discreteness of the ineffective subgroup).
\par
We can always reduce to the case of an almost
effective action of $\G$\,: at the level of \CR algebras,
this corresponds to substituting to $(\g,\q)$ its
\emph{effective quotient}, which is the
\CR algebra $(\g/ \mathfrak{a},\q/ \hat{\mathfrak{a}})$,
where $\mathfrak{a}$ is
the maximal ideal of $\g$ that is contained in
$\g_+$, and $\Hat{\mathfrak{a}}$ its
complexification in $\Hat{\g}$ (see \cite[\S 4]{MN:2005}).
\section{Finiteness and nondegeneracy conditions}
\label{sec:nondegeneracy}
We recall (see \cite[\S 13]{MN:2005}) that a $CR$ manifold $M$ is:
\begin{itemize}
\item of \emph{finite kind} (or \emph{finite type}) at $p\in M$ if
the higher order commutators of $C^\infty(M,HM)$, evaluated at
$p$, span $T_pM$;
\item \emph{holomorphically nondegenerate} at
$p\in M$ if there is no germ of nonzero holomorphic vector field
$X$ at $p$ that is tangent to $M$ at all points of a neighborhood
of $p$ (see \cite[\S 11.3]{BER:1999});
\item \emph{weakly nondegenerate} at
$p\in M$ if for each $\bar Z\in C^\infty(M,T^{0,1}M)$, with $\bar
Z_{|p}\neq 0$, there exist $m\in\mathbb N$ and $Z_1,\dots,Z_m\in
C^\infty(M,T^{1,0}M)$ such that:
\[
[Z_1,\dots,Z_m,\bar Z]_{|p}=[Z_1,[Z_2,\dots,[Z_m,\bar
Z]\dots]]_{|p} \not\in T^{1,0}_pM+T^{0,1}_pM\,;
\]
\item \emph{strictly nondegenerate} (or \emph{Levi nondegenerate})
at $p\in M$ if for each $\bar Z\in C^\infty(M,T^{0,1}M)$, with
$\bar Z_{|p}\neq 0$, there exists $Z'\in C^\infty(M,T^{1,0}M)$
such that:
\[
[Z',\bar Z]_{|p}\not\in T^{1,0}_pM+T^{0,1}_pM\,.
\]
\end{itemize}
\medskip\par
We now give analogous conditions for a \CR algebra. Let
$(\g,\q)$ be a \CR algebra of type $(n,k)$.
We say that $(\g,\q)$ is:
\begin{itemize}
\item
\emph{fundamental} if
${\HHH }_+$ generates $\g$ as a Lie algebra; or, equivalently,
$\q+\overline{\q}$ generates~$\hat\g$ as a complex
Lie algebra;
\item
\emph{ideal nondegenerate}, if there is no ideal $\mathfrak a$ of
$\g$ with ${\mathfrak a}\subset{\HHH }_+$ and ${\mathfrak
a}\not\subset\g_+$;
\item
\emph{weakly nondegenerate}, if
there is no  complex subalgebra
$\q'\neq\q$ of $\hat \g$
with $\q\subset\q'\subset\q+\overline{\q}$\, ;
\item
\emph{strictly nondegenerate}, if
$\left\{ X\in {\HHH }_+\,| \, [X,{\HHH }_+]\subset {\HHH }_+\right\}
\,=\, \g_+$.
\end{itemize}
Clearly we have the implications:
$$\boxed{
\text{{\sl strictly nondegenerate} $\Rightarrow$ {\sl weakly
nondegenerate} $\Rightarrow$ {\sl ideal nondegenerate.}}}
$$
Examples given in \cite[\S 2]{MN:2005} show that these conditions are
not equivalent.
\smallskip\par
For effective \CR algebras $(\g,\q)$, ideal nondegenerate
reduces to the statement:
$$\boxed{\text{
{\sl there are no nontrivial ideals of $\g$ contained in
${\HHH }_+$}}}$$
and this condition is actually equivalent to
\emph{effective and ideal nondegenerate}.
\\
The main motivation for considering \emph{ideal nondegenerate}
\CR algebras is the following (see \cite[\S 9]{MN:2005}):
\begin{thm}
Let $(\g,\q)$ be an ideal nondegenerate effective fundamental
\CR algebra.
Then $\g$ is finite dimensional.\qed
\end{thm}
\smallskip\par
We have the following (see \cite[\S 13]{MN:2005}):
\begin{thm}
Let $(M,HM,J)$ be a \CR manifold of type
$(n,k)$, homogeneous for the action of a Lie group $\G$, and
$(\g,\q)$ be the corresponding \CR algebra.
Then $(M,HM,J)$ is of finite kind (resp.\ holomorphically, weakly,
strictly nondegenerate) at all $p\in M$ if and only if the \CR algebra
$(\g,\q)$ is fundamental (resp.\ ideal, weakly, strictly
nondegenerate).\qed
\end{thm}
\section{Equivariant fibrations and reductions}\label{secequifib}
A \emph{morphism of \CR algebras}
$\phi:(\g,\q)\longrightarrow(\g',\q')$ is
a homomorphism $\phi:\g\longrightarrow\g'$ of real Lie algebras
whose complexification $\Hat{\phi}$ satisfies
$\Hat{\phi}(\q)\subset\q'$; it is a \emph{\CR submersion}
if $\phi(\g)+\g_+'=\g'$ and
$\Hat{\phi}(\q)+\q'\cap\bar{\q}'=\q'$.
\par
A \CR map $f:M_1\to M_2$ is a \emph{\CR fibration} if it is a smooth
submersion, i.e.\ $f_*(TM_1)=TM_2$, and a \emph{\CR submersion}, i.e.\
$f_*(HM_1)=HM_2$. In particular this is the case if $\q_1\subset\q_2$
and $f:M_1=\tilde{M}(\g,\q_1)\longrightarrow\tilde{M}(\g,\q_2)=M_2$ is
the natural projection.
\par
In this section we discuss morphisms of \CR algebras of the
special form
\begin{equation}\label{CRfibration}
\Phi_{\q}^{\q'}:
(\g,\q)\longrightarrow(\g,\q')\,,\quad\text{for $\q\subset\q'$.}
\end{equation}
They have been called in \cite{MN:2005}
\emph{$\g$-equivariant fibrations} and describe at the level
of \CR algebras the corresponding
$\G$-equivariant smooth \CR fibrations
$\tilde{M}(\g,\q)\longrightarrow\tilde{M}(\g,\q')$.
\smallskip\par
The \CR algebra $(\g'',\q'')=(\g'_+,\q\cap\bar\q')$
is called the \emph{fiber} of the $\g$-equivariant fibration
(\ref{CRfibration}).
It completely determines $\q'\supset\q$ by
$\q'=\q+\hat\g''$.
\\
The
$\g$-equivariant fibration
$\Phi_{\q}^{\q'}:(\g,\q)\rightarrow(\g,\q')$
will also be written as:
\begin{equation}\label{fibra}
\begin{CD}
(\g,\q) @>>(\g '',\q'')> (\g,\q')\, ,
\end{CD}
\end{equation}
with $(\g,\q)$ the \emph{total space}, $(\g,\q')$ \emph{the base},
and $(\g'',\q'')$ \emph{the fiber} of the fibration \eqref{fibra}.
\par
Now we describe two important examples of equivariant \CR fibrations
(see \cite[\S 5]{MN:2005}).
\subsection{Reduction to fundamental \CR algebras}
\begin{prop}\label{realred}
Every \CR algebra $(\g,\q)$ admits a unique
$\g$-equi\-va\-ri\-ant \CR fibration \eqref{fibra}
with a totally real base and a fundamental fiber.\qed
\end{prop}
We shall refer to the canonical $\g$-equivariant fibration of
Proposition \ref{realred} as the \emph{fundamental reduction} of
the \CR algebra $(\g,\q)$.\par Let $\hat{\mathfrak f}$ be the
subalgebra of $\hat\g$ generated by ${\q}+\bar\q$. Then the base
of the reduction is $(\g,\hat{\mathfrak f})$ and the fiber is
$(\g'',\q'')=(\g \cap\hat{\mathfrak f},\q)$.
\smallskip\par
Note that if
$\Phi:(\g,\q)\rightarrow(\g',\q')$
is a \CR submersion and $(\g,\q)$
is fundamental, then also $(\g',\q')$ is fundamental.
\par
For general \CR manifolds we have, by Frobenius Theorem:
\begin{prop}
Let $(M,HM,J)$ be a \CR manifold. Let $B\subset TM$ be the
subbundle generated by $HM$ and its commutators. If the rank of
$B$ is constant in a small neighborhood $U$ of a point $p\in M$,
then $B$ generates a smooth local foliation $U=\bigcup U_\alpha$
such that each leaf $U_\alpha$ is a \CR manifold of finite type
and its \CR dimension is the same as the \CR dimension of $M$.\par
Furthermore, if $M$ is
locally $\G$-homogeneous for a Lie group $\G$ of \CR
transformations, then every $U_\alpha$ is locally homogeneous and
the associated \CR algebra is $(\g'',\q'')$.\qed
\end{prop}
\medskip
\subsection{Reduction to weakly nondegenerate \CR algebras}
\begin{prop}\label{compred}
Every fundamental \CR algebra $(\g,\q)$ admits
a unique $\g$-equivariant \CR fibration \eqref{fibra}
with a weakly nondegenerate fundamental base and a totally complex
fiber.\qed
\end{prop}
We shall refer to the canonical $\g$-equivariant fibration of
Proposition \ref{compred} as the \emph{weakly nondegenerate
reduction} of the fundamental \CR algebra $(\g,\q)$.
\par
When $(\g,\q)$ is weakly degenerate, it was proved
in \cite{MN:2005} that there is a \CR-fibration
$\tilde{M}(\g,\q)\longrightarrow M'$ of the corresponding homogeneous
\CR manifold $\tilde{M}(\g,\q)$
on a \CR manifold $M'$ with
the same \CR codimension, having a non trivial complex fiber.
For homogeneous
simply connected \CR manifolds, the condition of
weak degeneracy
is in fact necessary and sufficient
for the existence of \CR fibrations with non trivial complex
fibers. Indeed, for general \CR manifolds, the existence
of a \CR fibration with non trivial complex fibers
implies weak degeneracy, as we have\,:
\begin{prop}
Let $M$ and $M'$ be \CR manifolds. We assume that
$M'$ is locally embeddable and that
there exists a \CR fibration
$M\stackrel{\pi}{\longrightarrow} M'$
with totally complex fibers of positive dimension.
Then $M$ is weakly degenerate.
\end{prop}
\begin{proof}
Let $f$ be any smooth \CR function defined on a neighborhood $U'$
of $p'\in M'$. Then $\pi^*f$ is a \CR function in
$U=\pi^{-1}(U')$, that is constant along the fibers of $\pi$.
Then, if $L\in\mathcal C^{\infty}(M,T^{1,0}M)$ is tangent to the
fibers of $\pi$ in $U$, we obtain that $[\bar Z_1\,,\,\,\hdots\,
,\,\,\bar Z_m,\,\,L]\left( \pi^*f\right)=0$ for every choice of
$\bar Z_1,\,\hdots ,\,\, \bar Z_m\in\mathcal
C^\infty(M,T^{0,1}M)$. Assume by contradiction that $M$ is weakly
nondegenerate at some $p$ with $\pi_{|p}=p'$. Then for some choice
of $\bar Z_1,\,\hdots ,\,\, \bar Z_m\in\mathcal
C^\infty(M,T^{0,1}M)$ we would have $v_{p}=[\bar
Z_1\,,\,\,\hdots\, ,\,\,\bar Z_m,\,\,L] \notin T_p^{1,0}M\oplus
T_p^{0,1}M$. Since the fibers of $\pi$ are totally complex,
$\pi_*(v_{p})\neq 0$. By the assumption that $M'$ is locally
embeddable at $p$, the real parts of the (locally defined) \CR
functions give local coordinates in $M'$ and therefore there is a
\CR function $f$ defined on a neighborhood $U'$ of $p'$ with
$v_p(\pi^*f)=\pi_*(v_{p})(f)\neq 0$. This gives a contradiction,
proving our statement.
\end{proof}
\subsection{The reduction diagram}
Let $(\g,\q)$ be a \CR algebra. Let us denote by
$\mathfrak f$ the Lie subalgebra of $\g$
generated by
${\HHH }_+(\g,\q)$ and by
$\q'$ the largest subalgebra of
$\hat\g$ such that
$\q\subset\q'\subset\q+\bar{\q}$.
Next we set $\g''=\f\cap\q'$
and $\q''=\q\cap\hat\g''$.
Then the canonical
fibrations discussed in this section can be
summarized in the diagram:
\begin{equation}
\begin{CD}
(\g,\q)
\\
@VV{\begin{CD}
(\f,\q) @>>(\g'',\q'')> (\f,\q')
\end{CD}}V
\\
(\g,\hat\f)
\end{CD}
\end{equation}
where $(\g,\hat{\mathfrak f})$ is totally real, $(\mathfrak f,\q)$ is
fundamental,
$(\g'',\q'')$ is totally complex,
$(\mathfrak f,\q')$ is fundamental and weakly nondegenerate.
\section{Minimal orbits in flag manifolds}
\label{sec:flag}
A \emph{complex flag manifold} is a coset space
$\X=\Hat{\mathbf{G}}/\mathbf{Q}$, where $\Hat{\mathbf{G}}$ is a
connected complex semisimple Lie group and $\mathbf{Q}$ is
parabolic in $\Hat{\mathbf{G}}$. The manifold $\X$ is a closed
complex projective variety which only depends on the Lie algebras
$\Hat{\g}$ of $\Hat{\mathbf{G}}$ and $\q$ of $\mathbf{Q}$\,: this
is a consequence of the fact that the center of a connected and
simply connected complex Lie group is contained in each of its
parabolic subgroups.
\par
A \emph{real form} of $\Hat{\mathbf{G}}$ is a real subgroup
$\mathbf G$ of $\Hat{\mathbf{G}}$ whose
Lie algebra $\g$ is a real form of $\hat\g$
(i.e.\ $\hat\g=\C\otimes_{\R} \g$).
The
real form $\g$ is the set of fixed points of
an anti-involution $\sigma$, called \emph{conjugation} in
$\Hat{\g}$\,:
$\g=
\mathrm{Fix}_{\hat\g}(\sigma)=\{X\in\hat\g\,|\,\sigma(X)=X\}$.
\par
A real form $\mathbf G$ acts on the complex flag manifold $\X$ by
left multiplication, and $\X$ decomposes into a disjoint union of
$\mathbf G$-orbits. In \cite{Wolf:1969} it is shown that there are
finitely many orbits including a unique one which is closed (hence
compact).
\par
This orbit $M$ has minimal dimension and is connected.
We will refer to $M$ as the \emph{minimal orbit} of $\G$ in $\X$.
In particular, the connected component of the identity
$\mathbf{G}^{\circ}$ of $\mathbf{G}$ is transitive on $M$. Thus,
while studying $M$, we can as well assume that
$\mathbf{G}=\mathbf{G}^{\circ}$ is connected.
\par
Moreover, up to conjugation, we can arrange  that the
closed orbit is
$M=\mathbf G\cdot o$, where $o=e\mathbf Q$.
We shall denote by $\mathbf G_+=\mathbf G\cap\mathbf Q$
the isotropy
subgroup of $\mathbf G$ at $o$ and
by $\g_+=\g\cap\q$ its Lie
algebra.
\par
The closed orbit $M$ has a natural $\G$-homogeneous \CR structure induced by its
embedding in the complex flag manifold $\X$ and the pair $(\g,\q)$
is the corresponding \CR algebra.
\section{Parabolic \CR algebras}
\label{sec:parabolic}
A \CR algebra $(\g,\q)$ is \emph{parabolic\/} if  $\g$ is finite
dimensional and $\q$ is a parabolic subalgebra of $\hat\g$. We
have (see \cite[\S 5]{AMN:2005}):
\begin{lem}
A parabolic \CR algebra $(\g,\q)$ is effective if
and only if $\g$ is semisimple and no simple ideal of $\hat\g$ is
contained in $\q\cap\bar{\q}$.\qed
\end{lem}
\noindent and
\begin{prop}
Let $(\g,\q)$ be an effective parabolic $CR$
algebra
and let $\g =\g_1\oplus \cdots \oplus\g_\ell$ be the
decomposition of $\g$ into the direct sum of its simple ideals.
Then:
\begin{enumerate}
\item $\q=\q_1\oplus \cdots \oplus\q_\ell$
where $\q_j=\q\cap\hat\g_j$ for $j=1,\hdots,\ell$;
\item for each $j=1,\hdots,\ell$, $(\g_j,\q_j)$
is an effective parabolic $CR$ algebra;
\item $(\g,\q)$ is fundamental (resp.\ ideal, weakly, strictly
nondegenerate) if and only if for each $j=1,\hdots,\ell$, the
$CR$ algebra $ (\g_j,\q_j)$ is fundamental (resp.\ ideal weakly,
strictly
nondegenerate).\qed
\end{enumerate}
\end{prop}
Thus in the following we can assume, with no loss of generality, that
$\g$ is simple and $\q$ is
parabolic in $\hat\g$.
\par\smallskip
To an effective parabolic \CR algebra $(\g,\q)$ we
associate a \CR manifold $M=M(\g,\q)$, unique
modulo isomorphisms,
defined as the orbit
$\G\cdot o$ in $\hat\G/\mathbf Q$, where:
\begin{itemize}
\item
$\hat\G$ is a connected and simply connected Lie
group with Lie algebra $\hat\g$\,;
\item
$\mathbf Q=\mathbf N_{\hat\G}(\q)$, the normalizer of $\q$ in $\hat\G$, is the
parabolic subgroup of $\hat\G$ with Lie algebra $\q$\,;
\item
$\G$ is the analytic real subgroup of $\hat\G$
with Lie algebra $\g$.
\end{itemize}
Note that for a parabolic $(\g,\q)$
the universal covering $\tilde{M}(\g,\q)$ of ${M}(\g,\q)$ admits a canonical
\CR structure, such that the covering map is a local \CR isomorphism.
If $\g
=\g_1\oplus \cdots \oplus\g_\ell$ is semisimple, then:
\begin{gather*}
M(\mathfrak g,\mathfrak q)\simeq M(\mathfrak g_1,\mathfrak q_1)
\times\cdots\times
M(\mathfrak g_\ell,\mathfrak q_\ell)\,;\\
\tilde{M}(\mathfrak g,\mathfrak q)\simeq
\tilde{M}(\mathfrak g_1,\mathfrak q_1)
\times\cdots\times
\tilde{M}(\mathfrak g_\ell,\mathfrak q_\ell).
\end{gather*}
\par\medskip
A real Lie subalgebra
$\mathfrak t$ of $\g$
is \emph{triangular} if all linear maps
$\mathrm{ad}_{\g}(X)\in\mathfrak{gl}_{\R}(\g)$ with $X\in\mathfrak
t$
can be simultaneously represented by triangular matrices in a suitable
basis of $\g$. All maximal triangular subalgebras of $\g$ are
conjugate by inner automorphisms (cf. \cite[\S 5.4]{Mostow:1961}).
A real Lie subalgebra of $\g$ containing a maximal triangular
subalgebra of $\g$ is called a $t$-\emph{subalgebra}.
\par
An effective parabolic \CR algebra $(\g,\q)$ will be
called \emph{minimal} if $\g_+={{\q}\cap{\g}}$ is
a $t$-subalgebra of $\g$.
\par
We observe that a maximal triangular subalgebra of $\g$ contains
a maximal Abelian subalgebra of semisimple elements having real
eigenvalues.
Hence:
\begin{prop}\label{propaf}
If $(\g,\q)$ is an effective parabolic minimal
\CR algebra, then there exists a Cartan subalgebra $\h$ of $\g$
with maximal vector part contained in $\q$.\qed
\end{prop}
A Cartan subalgebra with these properties is said to be \emph{adapted}
to $(\g,\q)$.
\par
The \emph{analytic Weyl group} of $\g$ with respect to $\h$ is
a subgroup of the Weyl group naturally isomorphic
to the quotient $\gr N_{\gr{Int}(\g)}(\h)/\gr Z_{\gr{Int}(\g)}(\h)$. Here $\gr Z$
denotes the centralizer, and $\gr{Int}(\g)$ is the group
of inner automorphisms of $\g$.
\begin{thm}\label{teoremaag}
Let $\g$ be a semisimple real Lie algebra and $\q$ a parabolic
subalgebra of its complexification $\Hat{\g}$. Then, up to
isomorphisms of \CR algebras, there is a unique parabolic minimal \CR algebra
$(\g',\q')$ with $\g'$ isomorphic to $\g$ and
$\q'$ isomorphic to $\q$.
\end{thm}
\begin{proof}
Fix a maximal triangular subalgebra $\mathfrak t$ of $\mathfrak
g$. Its complexification $\hat{\mathfrak t}$ is solvable,
therefore is contained in a maximal solvable subalgebra, i.e.\ a
Borel subalgebra $\mathfrak b$ of $\Hat{\mathfrak g}$. Modulo an
inner automorphism of $\Hat{\mathfrak g}$, we can assume that
$\mathfrak b\subset\mathfrak q$. The $CR$ algebra $(\mathfrak
g,\mathfrak q)$ is parabolic minimal. \par Let $\mathfrak q$,
$\mathfrak q'$ be parabolic subalgebras of $\Hat{\mathfrak g}$
such that $\mathfrak g_+={{\q}\cap{\g}}$ and $\mathfrak
g'_+={{\q'}\cap{\g}}$ are $t$-subalgebras of $\mathfrak g$. By an
inner automorphism of $\mathfrak g$, we can assume that $\mathfrak
g_+$ and  $\mathfrak g'_+$ contain the same maximal triangular
subalgebra $\mathfrak t$ of $\mathfrak g$, hence the same maximal
Abelian subalgebra of $\mathfrak g$ of semisimple elements having
real eigenvalues. Then, using another inner automorphism of
$\mathfrak g$, we can assume that $\mathfrak q$ and $\mathfrak q'$
contain the same maximal vectorial Cartan subalgebra $\mathfrak h$
of $\mathfrak g$. \par The inner automorphism of $\Hat{\mathfrak
g}$ transforming $\mathfrak q$ into $\mathfrak q'$ can now be
taken to be an element of the analytic Weyl group, leaving the
Cartan subalgebra $\mathfrak h$ and hence $\mathfrak g$ invariant.
It defines a $CR$ isomorphism between $(\mathfrak g,\mathfrak q)$
and $(\mathfrak g,\mathfrak q')$.
\end{proof}
Effective parabolic minimal \CR algebras correspond to minimal orbits.
In fact we have\,:
\begin{thm}
The \CR manifold $M(\g,\q)$
associated to an effective parabolic subalgebra
$(\g,\q)$
is compact if and only if $(\g, \q)$ is minimal.
\end{thm}
\begin{proof}
Indeed, since $\mathbf{G}$ is a linear group, a $\G$-homogeneous
space  $\G/\G_+$ is compact if and only if $\G_+$ contains a
maximal connected triangular subgroup (see \cite[Part II, Ch.5, \S
1.1]{GOV:1993}), i.e.\ if $\g_+$ is a $t$-subalgebra of $\g$.
\end{proof}
\section{Cross marked Satake diagrams}\label{secsatake}
We recall some facts about Satake diagrams and parabolic subgroups.
\par
Let $\g$ be a real semisimple Lie algebra and $\h$ a Cartan subalgebra
with maximal vector part. Fix a Cartan decomposition $\g=\la
k\oplus\la p$ such that
\[
\h=\h^+\oplus\h^-\text{ where }\h^+=\h\cap\la k\text{ and }
\h^-=\h\cap\la p
\]
and set $\h_\R=(i\h^+)\oplus\h^-$.
\par
The root system $\ro=\ro(\hat\g,\hat\h)$ of $\hat\g$ with respect
to $\hat\h$ is contained in $\h_\R^*$, the real dual of $\h_\R$. The
conjugation $\sigma$ of $\hat\g$ induces an involution on $\h_\R^*$
(that we still denote by $\sigma$), that preserves $\ro$.
\par
A root $\alpha\in\ro$ is called \emph{real} if
$\bar\alpha=\sigma(\alpha)=\alpha$,
\emph{imaginary} if $\bar\alpha=\sigma(\alpha)=-\alpha$.
We shall denote by $\ro_{\bullet}$ the set of imaginary roots
in $\ro$.
\par
For a Weyl chamber $C$, denote by $\ro^+(C)$ the set of positive
roots. If $C'$ is another Weyl chamber, let $w_{C,C'}$ be the
unique element of the Weyl group sending $C$ to $C'$. We have (see
\cite{Araki:1962}):
\begin{prop}\label{propcamera}
There exists a Weyl chamber $C$, unique modulo the action of the
analytic Weyl group, such that $\bar\alpha\succ 0$ for all
$\alpha\in\ro^+(C)\setminus\ro_\bullet$.\qed
\end{prop}
This Weyl chamber determines a basis $\B=\B(C)$ of $\ro$. The map
\begin{equation}\label{eps-C}
\epsilon_C=\sigma\circ w_{C,\bar C}
\end{equation}
is an involution of $\ro$ that preserves $C$, and thus $\B$.
Moreover, for every root $\alpha\in\B\setminus\ro_\bullet$,
there are integers $n_{\alpha,\beta}\geq 0$ such that
\begin{equation}\label{formulaconj}
\bar\alpha=\varepsilon_C(\alpha)\, +  \,
\sum_{\beta\in\B\cap\ro_{\bullet}}{n_{\alpha,\beta}
\beta}\,.
\end{equation}
\par
The \emph{Satake diagram} $\pmb{\mathcal{S}}$ of $\g$ is
obtained by the Dynkin diagram of $\hat\g$ whose nodes correspond to
the roots in $\B$ by painting black those corresponding to imaginary
roots and joining by a curved arrow those corresponding to
pairs of distinct nonimaginary roots $\alpha,\beta$ with
$\epsilon_C(\alpha)=\beta$.\par
There is a bijective correspondence between real semisimple Lie
algebras and Satake diagrams.
\par\smallskip
Let $\hat\g$ be a complex semisimple Lie algebra, $\hat\h$ any Cartan subalgebra,
$\ro$ the root system, $C$ a Weyl chamber,
$\ro^+$ and $\ro^-$ the corresponding sets of positive and negative roots,
$\B=\{\alpha_1,\ldots,\alpha_{\ell}\}$ the corresponding basis.
\par
Denote by $\hat\g^\alpha$ the eigenspace in $\hat\g$ for the root
$\alpha$.
If $\beta=\sum_i n_i\alpha_i$, $n_i\geq 0$, define
the \emph{support} of $\beta$:
$\supp(\beta)=\{\alpha_i\,|\,n_i>0\}$.
\par
Let $\Phi$ be a subset of $\B$. The set $\check\Phi^r$ of those
$\beta\in\ro^-$ for which $\supp(\beta)\cap\Phi=\emptyset$
is a closed system of roots
(i.e. $\beta_1,\beta_2\in\check\Phi^r$ and $\beta_1+\beta_2\in\ro$
$\Rightarrow$ $\beta_1+\beta_2\in\check\Phi^r$).
Then
\begin{equation}\label{eqparab}
\q_{\Phi}=
\hat\h\oplus\sum_{\beta\in{\ro}^+}{\Hat{\g}^{\beta}}\oplus
\sum_{\beta
\in\check\Phi^r}{\Hat{\g}^{\beta}}
\end{equation}
{is} a parabolic subalgebra of $\Hat{\g}$,
and every parabolic
subalgebra of $\Hat{\g}$ can be
described in this way.
More precisely:
\begin{prop}
Let $\hat\g$ be a complex semisimple Lie algebra and $\q$ a parabolic
subalgebra. Then there exist a Cartan subalgebra $\hat\h$, a Weyl
chamber $C$ and a subset $\Phi\subset\B$ such that $\q=\q_\Phi$.\qed
\end{prop}
To a parabolic subalgebra $\q_\Phi$ we associate a \emph{cross-marked
Dynkin diagram}, consisting of the Dynkin diagram of $\g$ with
cross-marks on the nodes corresponding to roots in $\Phi$.
The correspondence between isomorphism classes of parabolic subalgebras
and cross-marked Dynkin diagrams is bijective.
\par\smallskip
If $(\g,\q)$ is an effective parabolic minimal \CR algebra, these two
constructions are compatible, in the following sense (see \cite[\S 6]{AMN:2005}):
\begin{thm}
Let $(\g,\q)$ be an effective parabolic minimal \CR algebra. Then there
exist a Cartan subalgebra $\h$ of $\g$ with maximal vector part, a Weyl
chamber $C$ for $\ro(\hat\g,\hat\h)$ and a subset $\Phi$ of $\B(C)$
such that $\q=\q_\Phi$.
\end{thm}
Then to an effective parabolic minimal \CR algebra we can associate a
\emph{cross-marked Satake diagrams} $(\pmb{\mathcal{S}},\Phi)$ where
$\pmb{\mathcal{S}}$ is the Satake diagram of $\g$ and $\q=\q_\Phi$;
this correspondence is bijective, modulo isomorphisms.
\par
We set up some notations that we will use in the following sections.
If $\q=\q_\Phi$, let $\mathcal Q=\ro^+\cup\check\Phi^r$, $\mathcal
Q^r=\check\Phi^r\cup(-\check\Phi^r)$ and $\mathcal Q^n=\mathcal
Q\setminus\mathcal Q^r$.
To the partition $\mathcal Q= \mathcal Q^r\cup
\mathcal Q^n$ corresponds a direct sum decomposition
$\q=\q^r\oplus\q^n$, where
\[
\q^r=\hat\h\oplus\sum_{\beta\in\mathcal Q^r}\hat\g^\beta\quad
\text{and}\quad \q^n=\sum_{\beta\in\mathcal Q^n}\hat\g^\beta
\]
are the \emph{reductive} and \emph{nilpotent} part of $\q$.
\section{\CR fibrations for parabolic minimal \CR algebras}
\label{sec:crfibr}
In this section we discuss morphisms of $CR$
algebras of the special form $(\g,\q)\to(\g,\q')$, for
$\q\subset\q'$. They have been called in \cite{MN:2005}
\emph{$\g$-equivariant fibrations} and describe at the level of
$CR$ algebras the corresponding $\mathbf{G}$-equivariant smooth
fibrations $M(\g,\q) \to M(\g,\q')$.
\\
We keep the notation of the
previous sections. In particular, $\g$ is a semisimple real Lie
algebra, $\h$ a Cartan subalgebra of $\g$ with maximal vector
part, $\ro=\ro(\Hat{\g},\Hat{\mathfrak{h}})$, $C$ a Weyl chamber
adapted to the conjugation $\sigma$ in $\mathfrak{h}^*_{\R}$
induced by the real form $\g$ of $\Hat{\g}$
(as in Proposition \ref{propcamera}),
$\B=\B(C)$ is the set
of simple roots in $\ro^+=\ro^+(C)$.\par Let
$\Psi\subset\Phi\subset\B$. Then $\q_\Phi\subset\q_\Psi$ and the
identity on $\g$ defines a natural $\g$-equivariant morphism of
$CR$ algebras\,:
\begin{equation}\label{formulafiba}
\pi: (\g,\q_\Phi) \to (\g,\mathfrak q_\Psi)\, .
\end{equation}
Its fiber is the \CR algebra $(\g',\q')$ where\,:
\begin{equation}
\begin{cases}
\g'=\g\cap\q_{\Psi}=\g\cap\bar{\mathfrak{q}}_{\Psi}
\,, \quad
\Hat{\g}'=\q_{\Psi}\cap\bar{\q}_{\Psi}\\
\q'=\Hat{\q}_{\Phi}\cap\Hat{\g}'=
\Hat{\q}_{\Phi}\cap\q_{\Psi}\cap\bar{\q}_{ \Psi}=
\q_{\Phi}\cap\bar{\q}_{\Psi}\, .
\end{cases}
\end{equation}
Denote by $\ro'$ and $\Q'$ the sets of roots $\alpha\in\ro$ for
which  $\Hat{\g}^{\alpha}$ is contained in $\hat\g'$ and
$\hat\q'$, respectively\,:
\begin{equation}
\begin{cases}
\ro'=\Q_\Psi\cap\bar{\Q}_\Psi \\
\Q'=\Q_\Phi\cap\bar\Q_\Psi\, .
\end{cases}
\end{equation}
\par
The \CR algebra $(\g',\q')$ in general is neither parabolic nor
effective. However its effective quotient \emph{is} parabolic
minimal. Define\,:
\begin{equation}
\begin{cases}
\ro''=\ro'\cap(-\ro')\,=\,\Q^r_{\Psi}\cap
\bar{\Q}^r_{\Psi}\\
\Q''=\Q'\cap\ro''
\end{cases}
\end{equation}
and set\,:
\begin{equation}
\begin{cases}
\hat\g''=\hat\h\oplus\bigoplus_{\alpha\in\mathcal
R''}\hat\g^\alpha\\
\q''=\q'\cap \hat\g''
\end{cases}
\end{equation}
Then $\ro''$ is $\sigma$-invariant, $\hat\g''=\mathfrak
q_\Psi^r\cap\bar\q_\Psi^r$ is reductive and $\mathfrak q''$ is
parabolic in $\hat\g''$.
Furthermore $\B''=\B\cap\ro''$ is a basis of $\ro''$.
\par
Hence we obtain (see \cite[\S 7]{AMN:2005}):
\begin{prop}\label{propgef}
The \CR algebra $(\g'',\q'')$ is parabolic minimal. Its
cross-marked Satake diagram $({\pmb{\mathcal S}''},\Phi'')$ is the
subdiagram of $(\pmb{\mathcal{S}},\Phi)$ consisting of the simple
roots $\alpha$ such that either one of the following conditions
holds\,:
\begin{enumerate}
\item $\alpha\in\ro_{\bullet}\setminus\Psi$; \item
$\alpha\not\in\ro_\bullet$ and
$\left(\{\alpha\}\cup\supp(\bar\alpha)
\right)\cap\Psi=\emptyset$.
\end{enumerate}\par
The cross-marks are left on the nodes corresponding to roots in
$\Phi\cap\B''$.\qed
\end{prop}
\par
We say that a Satake diagram is \emph{$\sigma$-connected} if
either it is connected or consists of two connected components,
joined by curved arrows.
\begin{thm}\label{thm:gef}
Let \eqref{formulafiba} be a $\g$-equivariant fibration. Then the
effective quotient of its fiber is the parabolic minimal \CR
algebra whose cross-marked Satake diagram consists of the union of
all $\sigma$-connected components of the diagram ${\pmb{\mathcal
S}''}$ described in Proposition \ref{propgef}, containing at least
one cross-marked node.\qed
\end{thm}
\par
Recall that a $\g$-equivariant morphism of $CR$ algebras
\eqref{formulafiba} is a \CR-fibration if the quotient map $
\q_{\Phi}/\left(\q_{\Phi}\cap\bar{\q} _{\Phi}\right) \to
\q_{\Psi}/\left(\q_{\Psi}\cap\bar{\q} _{\Psi}\right) $ is
onto.
\par
The condition that \eqref{formulafiba} is a \CR-fibration is
equivalent to the fact that the quotient map $M(\g,\q_{\Phi})\to
M(\g,\q_{\Psi})$ is a \CR fibration with fiber $M(\g'',\q'')$.
This condition is always satisfied if $(\g,\q_{\Psi})$ is totally
real, indeed in this case
$\q_{\Psi}/\left(\q_{\Psi}\cap\overline{\q} _{\Psi}\right)=0$.
\par
\begin{exam}
Let $\g=\mathfrak{su}(1,3)$ and let
$\Phi=\{\alpha_1,\alpha_2\}$, $\Psi=\{\alpha_1\}$. Then
the cross-marked Satake diagrams corresponding to the $CR$ algebra
$(\g,\q_\Phi)$, the basis $(\g,\q_\Psi)$
and the corresponding effective fiber are given by:
\par\bigskip
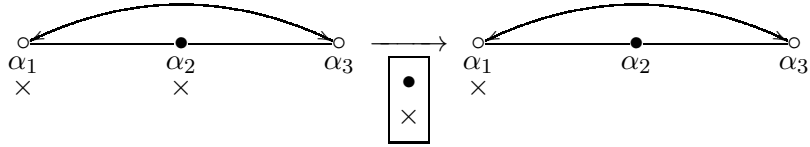
\[\begin{CD}
{\xymatrix@M=0pt @R=2pt{
\circ\ar @{-} [rr] \ar @/^15pt/@{<->}[rrrr]&&\bullet\ar@{-}[rr]&&\circ\\
\alpha_1&&\alpha_2&&\alpha_3\\
\times&&\times}}@>>{\boxed{\begin{matrix}
\bullet\\
\times
\end{matrix}}}>
{\xymatrix@M=0pt @R=2pt{
\circ\ar @{-} [rr] \ar @/^15pt/@{<->}[rrrr]&&\bullet\ar@{-}[rr]&&\circ\\
\alpha_1&&\alpha_2&&\alpha_3\\
\times}}\end{CD}
\]
In the case $\Psi=\{\alpha_2\}$ we have instead:
\par\bigskip
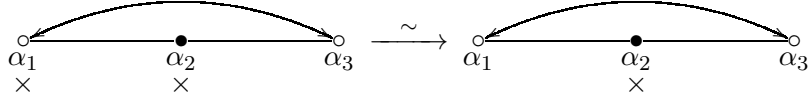
\[\begin{CD}
{\xymatrix@M=0pt @R=2pt{
\circ\ar @{-} [rr] \ar @/^15pt/@{<->}[rrrr]&&\bullet\ar@{-}[rr]&&\circ\\
\alpha_1&&\alpha_2&&\alpha_3\\
\times&&\times}}@>{\sim}>>
{\xymatrix@M=0pt @R=2pt{
\circ\ar @{-} [rr] \ar @/^15pt/@{<->}[rrrr]&&\bullet\ar@{-}[rr]&&\circ\\
\alpha_1&&\alpha_2&&\alpha_3\\
&&\times}}\end{CD}
\]
The fiber is trivial and the map is a $CR$ morphism, but not a
$CR$ isomorphism. The corresponding map
$M(\g,\q_{\Phi})\to M(\g,\q_{\Psi})$ is
an analytic diffeomorphism and a $CR$ map, but not a $CR$
diffeomorphism.
\end{exam}
\section{The fundamental reduction for parabolic minimal \CR algebras}
\label{sec:fundamental}
We give a criterion to read off the property of being fundamental
from the cross-marked Satake diagram\,:
\begin{prop}\label{propea}
An effective parabolic minimal \CR algebra $(\g,\q_{\Phi})$ is
fundamental if and only if its corresponding cross-marked Satake
diagram $(\pmb{\mathcal S},\Phi)$ has the property:
\[
\alpha\in\Phi\setminus\ro_{\bullet}
 \Longrightarrow \epsilon_C(\alpha)\notin \Phi\, .
\]
\end{prop}
Here $\epsilon_C$ is the involution in $\mathcal B(C)$ defined in
(\ref{eps-C}).
\begin{proof}
Assume that $\alpha_1$ and $\alpha_2=\epsilon_C(\alpha_1)$ both
belong to $\Phi$, and let $\Psi=\{\alpha_1,\alpha_2\}$. Then
$\Psi\subset\Phi$ and hence $\q_{\Phi}\subset\q_{\Psi}$. To show
that $(\g,\q_{\Phi})$ is not fundamental, it is sufficient to
check that $\q_{\Psi}=\overline{\q}_{\Psi}$. To this aim it
suffices to verify that $\mathcal Q^n_{\Psi}=\overline{\mathcal
Q}^{\, n}_{\Psi}$. Let $\mathcal
B(C)=\{\alpha_1,\alpha_2,\hdots,\alpha_\ell\}$. Every
root $\alpha\in\mathcal Q^n_{\Psi}$ can be written in the form
$\alpha=\sum_{i=1}^\ell{k_i\alpha_i}$ with $k_1+k_2>0$. Since $C$
is adapted to the conjugation $\sigma$, using \eqref{formulaconj}
we obtain\,:
\[
\bar\alpha=\sum_{i=1}^\ell{k_i\epsilon_C(\alpha_i)}+
\sum_{\beta\in\mathcal B\cap\mathcal R_\bullet}{k_{\alpha,\beta}\beta}
=\sum_{i=1}^\ell{k'_i\alpha_i},
\]
with $k'_1+k'_2=k_2+k_1>0$, showing that also
$\bar\alpha\in\mathcal Q^n_{\Psi}$. This shows that the condition
is necessary.
\par
Assume vice versa that there exists a proper parabolic subalgebra
$\q'$ of $\hat\g$ with $\q_{\Phi}\subset\q'=\overline{\q}\, '$.
Then $\q'=\q_{\Psi}$ for some $\Psi\subset\Phi$,
$\Psi\neq\emptyset$. Since $\overline{\mathcal Q}^{\,
n}_{\Psi}=\mathcal Q^n_{\Psi}\subset\mathcal R^+(C)$, we have
$\Psi\cap\ro_{\bullet}=\emptyset$. Hence, again by
\eqref{formulaconj}, we obtain that $\epsilon_C(\alpha)\in\Psi$
for all $\alpha\in\Psi$.
\end{proof}
\par
From Propositions \ref{propea} and \ref{propgef} we obtain\,:
\begin{thm}
Let $(\g,\q_{\Phi})$ be an effective parabolic minimal \CR algebra
and let $(\pmb{\mathcal S},\Phi)$ be its corresponding
cross-marked Satake diagram. Let
\[
\Psi=\{\alpha\in\Phi\setminus\mathcal{R}_{\bullet}\, | \,
\epsilon_C(\alpha)\in\Phi\}.
\]
Then
\begin{enumerate}
\item The diagram $\pmb{\mathcal S}'$ obtained from $\pmb{\mathcal
S}$ by erasing all the nodes corresponding to the roots in $\Psi$
and the lines and arrows issued from them is still a Satake
diagram, corresponding to a semisimple real Lie algebra $\g'$.
\item $(\g,\q_{\Psi})$ is a totally real effective parabolic
minimal \CR algebra. \item The natural map $(\g,\q_{\Phi})\rightarrow
(\g,\q_{\Psi})$, defined by the inclusion
$\q_{\Phi}\subset\q_{\Psi}$, is a $\mathfrak g$-equivariant \CR
fibration. The effective quotient of its fiber is the fundamental
parabolic minimal \CR algebra $(\g'',\q_{\Phi'})$, associated to
the cross-marked Satake diagram $(\pmb{\mathcal{S}}'',\Phi')$,
where $\Phi'=\Phi\setminus\Psi$ and $\pmb{\mathcal{S}}''$ is the
union of the $\sigma$-connected components of $\pmb{\mathcal{S}}'$
that contain some root of $\Phi'$.\qed
\end{enumerate}\end{thm}
\begin{exam}
Let $\g\simeq\mathfrak{su}(2,2)$ and let
$\Phi=\{\alpha_2,\alpha_3\}$ (we refer to the diagram below).
We have $\epsilon_C(\alpha_i)=\alpha_{4-i}$ for $i=1,2,3$
and hence
$\Psi=\{\alpha\in\Phi\, | \,
\epsilon_C(\alpha)\in\Phi\}=\{\alpha_2\}$.
In particular $(\g,\q_{\{\alpha_2,\alpha_3\}})$
is not fundamental. We obtain
a $\g$-equivariant
$CR$ fibration
$(\g,\q_{\{\alpha_2,\alpha_3\}})\to
(\g,\q_{\{\alpha_2\}})$ with
fundamental fiber
$(\g',\q'_{\{\alpha_3\}})$, with
$\g'\simeq\mathfrak{sl}(2,\C)$.
\[\begin{CD}
{\xymatrix@M=0pt @R=2pt{
\circ\ar @{-}[rr]\ar @/^15pt/@{<->}[rrrr]&&\circ\ar @{-}[rr]&&\circ\\
\alpha_1&&\alpha_2&&\alpha_3\\
&&\times&&\times}}@>{\begin{matrix}
\boxed{\xymatrix@M=0pt @R=2pt{
\quad\\
\quad \\
&\times\\
\circ\ar @/^10pt/@{<->}[r]&\circ\\
\alpha_1&\alpha_3
}}
\quad\\
\\
\end{matrix}}>>
 {\xymatrix@M=0pt @R=2pt{
\circ\ar @{-}[rr]\ar @/^15pt/@{<->}[rrrr]&&\circ\ar @{-}[rr]&&\circ\\
\alpha_1&&\alpha_2&&\alpha_3\\
&&\times}}\end{CD}
\]
\end{exam}
\def\MR#1{\href{http://www.ams.org/mathscinet-getitem?mr=#1}{\mbox{MR\;#1}}}
\providecommand{\bysame}{\leavevmode\hbox to3em{\hrulefill}\thinspace}
\providecommand{\href}[2]{#2}

\end{document}